\documentclass [10pt,reqno]{amsart}

\usepackage {hannah}
\usepackage{pstricks}
\usepackage{pst-plot}
\usepackage {color,graphicx}
\usepackage {multirow}
\usepackage{texdraw}
\usepackage[all]{xy}
\newcommand{\tom}[1]{}
\newcommand{\lang}[1]{}
\newcommand{\kurz}[1]{#1}

\DeclareMathOperator {\Int}{Int}

\DeclareMathOperator {\val}{val}
\DeclareMathOperator {\trop}{trop}
\DeclareMathOperator {\Trop}{Trop}

\DeclareMathOperator {\newton}{N}

\newcommand{\ca}{{\mathcal{A}}}
\CompileMatrices     % Diagramme sollen nur nach Aenderungen neu compiliert werden

\newcommand{\Pa}{
  \begin{texdraw}
    \drawdim cm \relunitscale 0.4
    \move (0 0) \lvec (4 0) \lvec (0 2) \lvec (0 0)
    \move (0 0) \fcir f:0 r:0.1
    \move (1 0) \fcir f:0 r:0.1
    \move (2 0) \fcir f:0 r:0.1
    \move (3 0) \fcir f:0 r:0.1
    \move (4 0) \fcir f:0 r:0.1
    \move (0 1) \fcir f:0 r:0.1
    \move (1 1) \fcir f:0 r:0.1
    \move (2 1) \fcir f:0 r:0.1
    \move (0 2) \fcir f:0 r:0.1
  \end{texdraw}
}
\newcommand{\Pb}{
  \begin{texdraw}
    \drawdim cm \relunitscale 0.4
    \move (0 0) \lvec (2 0) \lvec (2 2) \lvec (0 2) \lvec (0 0)
    \move (0 0) \fcir f:0 r:0.1
    \move (1 0) \fcir f:0 r:0.1
    \move (2 0) \fcir f:0 r:0.1
    \move (0 1) \fcir f:0 r:0.1
    \move (1 1) \fcir f:0 r:0.1
    \move (2 1) \fcir f:0 r:0.1
    \move (0 2) \fcir f:0 r:0.1
    \move (1 2) \fcir f:0 r:0.1
    \move (2 2) \fcir f:0 r:0.1
  \end{texdraw}
}
\newcommand{\Pc}{
  \begin{texdraw}
    \drawdim cm \relunitscale 0.4
    \move (0 0) \lvec (3 0) \lvec (0 3) \lvec (0 0)
    \move (0 0) \fcir f:0 r:0.1
    \move (1 0) \fcir f:0 r:0.1
    \move (2 0) \fcir f:0 r:0.1
    \move (3 0) \fcir f:0 r:0.1
    \move (0 1) \fcir f:0 r:0.1
    \move (1 1) \fcir f:0 r:0.1
    \move (2 1) \fcir f:0 r:0.1
    \move (0 2) \fcir f:0 r:0.1
    \move (1 2) \fcir f:0 r:0.1
    \move (0 3) \fcir f:0 r:0.1
  \end{texdraw}
}
\newcommand{\Pca}{
  \begin{texdraw}
    \drawdim cm \relunitscale 0.4
    \move (0 0) \lvec (3 0) \lvec (1 2) \lvec (0 2) \lvec (0 0)
    \move (0 0) \fcir f:0 r:0.1
    \move (1 0) \fcir f:0 r:0.1
    \move (2 0) \fcir f:0 r:0.1
    \move (3 0) \fcir f:0 r:0.1
    \move (0 1) \fcir f:0 r:0.1
    \move (1 1) \fcir f:0 r:0.1
    \move (2 1) \fcir f:0 r:0.1
    \move (0 2) \fcir f:0 r:0.1
    \move (1 2) \fcir f:0 r:0.1
  \end{texdraw}
}
\newcommand{\Pcb}{
  \begin{texdraw}
    \drawdim cm \relunitscale 0.4
    \move (0 0) \lvec (2 0) \lvec (2 1) \lvec (1 2) \lvec (0 2) \lvec (0 0)
    \move (0 0) \fcir f:0 r:0.1
    \move (1 0) \fcir f:0 r:0.1
    \move (2 0) \fcir f:0 r:0.1
    \move (0 1) \fcir f:0 r:0.1
    \move (1 1) \fcir f:0 r:0.1
    \move (2 1) \fcir f:0 r:0.1
    \move (0 2) \fcir f:0 r:0.1
    \move (1 2) \fcir f:0 r:0.1
  \end{texdraw}
}
\newcommand{\Pcc}{
  \begin{texdraw}
    \drawdim cm \relunitscale 0.4
    \move (0 0) \lvec (3 0) \lvec (1 2) \lvec (0 1) \lvec (0 0)
    \move (0 0) \fcir f:0 r:0.1
    \move (1 0) \fcir f:0 r:0.1
    \move (2 0) \fcir f:0 r:0.1
    \move (3 0) \fcir f:0 r:0.1
    \move (0 1) \fcir f:0 r:0.1
    \move (1 1) \fcir f:0 r:0.1
    \move (2 1) \fcir f:0 r:0.1
    \move (1 2) \fcir f:0 r:0.1
  \end{texdraw}
}
\newcommand{\Pcd}{
  \begin{texdraw}
    \drawdim cm \relunitscale 0.4
    \move (1 0) \lvec (2 0) \lvec (2 1) \lvec (1 2) 
    \lvec (0 2) \lvec (0 1) \lvec (1 0)
    \move (1 0) \fcir f:0 r:0.1
    \move (2 0) \fcir f:0 r:0.1
    \move (0 1) \fcir f:0 r:0.1
    \move (1 1) \fcir f:0 r:0.1
    \move (2 1) \fcir f:0 r:0.1
    \move (0 2) \fcir f:0 r:0.1
    \move (1 2) \fcir f:0 r:0.1
  \end{texdraw}
}
\newcommand{\Pce}{
  \begin{texdraw}
    \drawdim cm \relunitscale 0.4
    \move (0 0) \lvec (2 0) \lvec (2 1) \lvec (1 2) 
    \lvec (0 1) \lvec (0 0)
    \move (0 0) \fcir f:0 r:0.1
    \move (1 0) \fcir f:0 r:0.1
    \move (2 0) \fcir f:0 r:0.1
    \move (0 1) \fcir f:0 r:0.1
    \move (1 1) \fcir f:0 r:0.1
    \move (2 1) \fcir f:0 r:0.1
    \move (1 2) \fcir f:0 r:0.1
  \end{texdraw}
}
\newcommand{\Pcf}{
  \begin{texdraw}
    \drawdim cm \relunitscale 0.4
    \move (0 0) \lvec (2 0) \lvec (1 2) \lvec (0 2) \lvec (0 0)
    \move (0 0) \fcir f:0 r:0.1
    \move (1 0) \fcir f:0 r:0.1
    \move (2 0) \fcir f:0 r:0.1
    \move (0 1) \fcir f:0 r:0.1
    \move (1 1) \fcir f:0 r:0.1
    \move (0 2) \fcir f:0 r:0.1
    \move (1 2) \fcir f:0 r:0.1
  \end{texdraw}
}
\newcommand{\Pcg}{
  \begin{texdraw}
    \drawdim cm \relunitscale 0.4
    \move (0 0) \lvec (3 0) \lvec (1 2) \lvec (0 0)
    \move (0 0) \fcir f:0 r:0.1
    \move (1 0) \fcir f:0 r:0.1
    \move (2 0) \fcir f:0 r:0.1
    \move (3 0) \fcir f:0 r:0.1
    \move (1 1) \fcir f:0 r:0.1
    \move (2 1) \fcir f:0 r:0.1
    \move (1 2) \fcir f:0 r:0.1
  \end{texdraw}
}
\newcommand{\Pch}{
  \begin{texdraw}
    \drawdim cm \relunitscale 0.4
    \move (0 0) \lvec (2 0) \lvec (2 1) \lvec (1 2) \lvec (0 0)
    \move (0 0) \fcir f:0 r:0.1
    \move (1 0) \fcir f:0 r:0.1
    \move (2 0) \fcir f:0 r:0.1
    \move (1 1) \fcir f:0 r:0.1
    \move (2 1) \fcir f:0 r:0.1
    \move (1 2) \fcir f:0 r:0.1
  \end{texdraw}
}
\newcommand{\Pci}{
  \begin{texdraw}
    \drawdim cm \relunitscale 0.4
    \move (1 0) \lvec (2 0) \lvec (2 1) \lvec (1 2) 
    \lvec (0 1) \lvec (1 0)
    \move (1 0) \fcir f:0 r:0.1
    \move (2 0) \fcir f:0 r:0.1
    \move (0 1) \fcir f:0 r:0.1
    \move (1 1) \fcir f:0 r:0.1
    \move (2 1) \fcir f:0 r:0.1
    \move (1 2) \fcir f:0 r:0.1
  \end{texdraw}
}
\newcommand{\Pcj}{
  \begin{texdraw}
    \drawdim cm \relunitscale 0.4
    \move (1 0) \lvec (2 1) \lvec (1 2) 
    \lvec (0 1) \lvec (1 0)
    \move (1 0) \fcir f:0 r:0.1
    \move (0 1) \fcir f:0 r:0.1
    \move (1 1) \fcir f:0 r:0.1
    \move (2 1) \fcir f:0 r:0.1
    \move (1 2) \fcir f:0 r:0.1
  \end{texdraw}
}
\newcommand{\Pck}{
  \begin{texdraw}
    \drawdim cm \relunitscale 0.4
    \move (1 0) \lvec (2 0) \lvec (1 2) \lvec (0 1) \lvec (1 0)
    \move (1 0) \fcir f:0 r:0.1
    \move (2 0) \fcir f:0 r:0.1
    \move (0 1) \fcir f:0 r:0.1
    \move (1 1) \fcir f:0 r:0.1
    \move (1 2) \fcir f:0 r:0.1
  \end{texdraw}
}
\newcommand{\Pcl}{
  \begin{texdraw}
    \drawdim cm \relunitscale 0.4
    \move (0 0) \lvec (2 0) \lvec (1 2) \lvec (0 0) 
    \move (0 0) \fcir f:0 r:0.1
    \move (1 0) \fcir f:0 r:0.1
    \move (2 0) \fcir f:0 r:0.1
    \move (1 1) \fcir f:0 r:0.1
    \move (1 2) \fcir f:0 r:0.1
  \end{texdraw}
}
\newcommand{\Pcm}{
  \begin{texdraw}
    \drawdim cm \relunitscale 0.4
    \move (2 0) \lvec (1 2) \lvec (0 1) \lvec (2 0) 
    \move (2 0) \fcir f:0 r:0.1
    \move (1 1) \fcir f:0 r:0.1
    \move (1 2) \fcir f:0 r:0.1
    \move (0 1) \fcir f:0 r:0.1
  \end{texdraw}
}

\title [The tropical $j$-invariant]{The tropical $j$-invariant}

\author {Eric Katz, Hannah Markwig, Thomas Markwig}
\address {Eric Katz, Department of Mathematics, The University of Texas
  at Austin, 1 University Station, C1200, Austin, TX 78712}
\email {eekatz@math.utexas.edu}
\address{Hannah Markwig, University of Michigan, Department of Mathematics, 2074 East Hall, 530 Church Street, Ann Arbor, MI 48109-1043}
\email{markwig@umich.edu}
\address {Thomas Markwig, Fachbereich Mathematik, Technische Universit\"at Kaiserslautern, Postfach 
3049, 67653 Kaiserslautern,
Germany}
\email {keilen@mathematik.uni-kl.de}
\thanks {\emph {2000 Mathematics Subject Classification:} 14H52, 51M20}
\thanks{The third author would like to thank the Institute for Mathematics and its Applications in 
Minneapolis for hospitality.}

\begin {document}

  \begin {abstract}  
    If $(Q,\ca)$ is a marked polygon with one interior point, then a
    general polynomial $f\in\K[x,y]$ with support $\ca$ defines an
    elliptic curve $C_f$ on the toric surface $X_\ca$. If $\K$ has a
    non-archimedean valuation into $\R$ we can tropicalize $C_f$ to get a
    tropical curve $\Trop(C_f)$. If the Newton subdivision induced by
    $f$ is a triangulation, then $\Trop(C_f)$ will be a graph of genus
    one and we show that the lattice length of the cycle of that graph
    is the negative of the valuation of the $j$-invariant of $C_f$. 
  \end {abstract}

  \maketitle

  \section{Introduction}

  Previous work by Grisha Mikhalkin \cite{Mik06}, by Michael Kerber
  and Hannah Markwig 
  \cite{KM08} and by Magnus Vigeland \cite{Vig04} shows that the length of
  the cycle of a tropical curve of genus one has properties which one
  classically attributes to the $j$-invariant of an elliptic
  curve without giving a direct link between these two numbers. 
  In \cite{KMM08} we established such a direct link for plane cubics
  by showing that the tropicalization of the $j$-invariant is \emph{in
    general} the negative of the cycle length. In the present paper we
  generalize this result to elliptic curves on other toric surfaces
  using the same methods.

  More precisely, if $(Q,\ca)$ is a marked polygon with one interior point, then a
  general polynomial $f\in\K[x,y]$ with support $\ca$ defines an
  elliptic curve $C_f$ on the toric surface $X_\ca$. If $\K$ has a
  non-archimedean valuation we can tropicalize $C_f$ to get a
  tropical curve $\Trop(C_f)$. If the Newton subdivision induced by
  $f$ is a triangulation, then $\Trop(C_f)$ will be a graph of genus
  one and we show in our main result in Theorem \ref{thm:main} that
  the lattice length of the cycle of the graph 
  is the negative of the valuation of the $j$-invariant. 

  In the case where the triangulation is \emph{unimodular}, i.e.\ all
  the triangles have area $\frac{1}{2}$, this result was
  independently derived by David Speyer \cite[Proposition 9.2]{Spe07} using  Tate
  uniformization of elliptic curves. David Speyer's result is more
  general though in the sense that it applies to curves in arbitrary toric
  varieties.

   This paper is organized as follows. In Section
   \ref{sec:toricsurfaces} we consider toric surfaces defined by a marked
   lattice polygon with one interior point, we recall the
   classification of these polygons and we consider the impact on the
   $j$-invariant for the corresponding elliptic curves. Section
   \ref{sec:tropicalization} recalls the notion of tropicalization and
   of plane tropical curves. We then introduce in Section
   \ref{sec:jinv} the notion of tropical $j$-invariant and give a
   formula to compute it. Section \ref{sec:unimodular} shows that the
   tropical $j$-invariant is preserved by integral unimodular affine
   transformations. With this preparation we are able to state our
   main result in Section \ref{sec:main}. Section \ref{sec:reduction}
   is then devoted to the reduction of the proof to considering only
   three marked polygons and Section \ref{sec:cases} shows how these
   three cases can be dealt with using procedures from the
   \textsc{Singular} library \texttt{jinvariant.lib} (see
   \cite{KMM07})  which is available via the URL 
   \begin{center}
     http://www.mathematik.uni-kl.de/\textasciitilde keilen/en/jinvariant.html.
   \end{center}
   The actual computations are done using \texttt{polymake}
   \cite{GJ97}, TOPCOM \cite{Ram02} and \textsc{Singular} \cite{GPS05}.
   The tropical curves in this paper and their Newton subdivisions were
   produced using the 
   procedure \texttt{drawtropicalcurve} from the 
   \textsc{Singular} library \texttt{tropical.lib} (see \cite{JMM07a}) which
   can be obtained via the URL
   \begin{center}
     http://www.mathematik.uni-kl.de/\textasciitilde keilen/en/tropical.html.
   \end{center}
   
   The authors would like to thank Vladimir Berkovich, Jordan
   Ellenberg, Bjorn Poonen, David Speyer, Charles Staats, Bernd Sturmfels
   and John Voigt for valuable discussions. 

%%%%%%%%%%%%%%%%%%%%%%%%%%%%%%%%%%%%%%%%%%%%%%%%%%%%%%%%%%%%%%%%%%%%%%%%%

   \section{Toric surfaces}\label{sec:toricsurfaces}
   
   Throughout this paper we consider mainly marked 
   polygons $(Q,\ca)$ such that $Q$ contains a single
   interior lattice point, where by a \emph{marked polygon} we mean a
   convex lattice polygon $Q$ in $\R^2$  together with a subset
   $\ca\subseteq Q\cap\Z^2$ of the lattice points of $Q$ containing at
   least the vertices of $Q$ (cf.~\cite[Section~2.A]{GKZ94}).
   Fixing a base field $\K$ 
   such a polygon defines a polarized \emph{toric surface}
   \begin{displaymath}
     X_\ca\subset\P_\K^{|\ca|-1}.
   \end{displaymath}
   In the torus
   $(\K^*)^2\subset X_\ca$ the hyperplane section, say $C_f$, defined by the
   linear form $\sum_{(i,j)\in\ca} a_{ij}\cdot z_{ij}$ is the
   vanishing locus of the Laurent polynomial
   \begin{displaymath}
     f=\sum_{(i,j)\in\ca}a_{ij}\cdot x^iy^j
   \end{displaymath}
   (cf.~\cite[Chapter~5]{GKZ94}). Since the
   arithmetical genus of the 
   hyperplane sections is the number of interior lattice points of
   $Q$ (cf.~\cite[p.~91]{Ful93}),
   the general hyperplane section will be a smooth \emph{elliptic
     curve}. The \emph{$j$-invariant} of such a curve is an element of
   the base field which characterizes the curve up to isomorphism.

\begin{figure}[h]
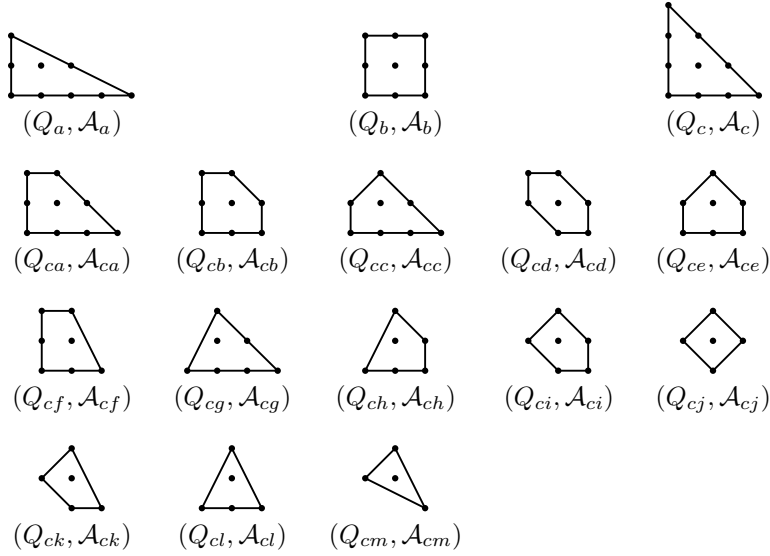

     \centering
     \begin{tabular}{c@{\hspace*{0.5cm}}c@{\hspace*{0.5cm}}c@{\hspace*{0.5cm}}c@{\hspace*{0.5cm}}c}
       \Pa && \Pb && \Pc  \\
       $(Q_a,\ca_a)$ && $(Q_b,\ca_b)$ && $(Q_c,\ca_c)$ \\[0.4cm]
       \Pca & \Pcb & \Pcc & \Pcd & \Pce\\
       $(Q_{ca},\ca_{ca})$ &$(Q_{cb},\ca_{cb})$ &$(Q_{cc},\ca_{cc})$ &$(Q_{cd},\ca_{cd})$ &$(Q_{ce},\ca_{ce})$ \\[0.4cm]
       \Pcf & \Pcg & \Pch & \Pci & \Pcj\\
       $(Q_{cf},\ca_{cf})$ &$(Q_{cg},\ca_{cg})$ &$(Q_{ch},\ca_{ch})$ &$(Q_{ci},\ca_{ci})$ &$(Q_{cj},\ca_{cj})$ \\[0.4cm]
       \Pck & \Pcl & \Pcm \\
       $(Q_{ck},\ca_{ck})$ &$(Q_{cl},\ca_{cl})$ &$(Q_{cm},\ca_{cm})$ \\[0.4cm]
     \end{tabular}     
     \caption{The $16$ convex lattice polygons with one interior
       lattice point}
     \label{fig:polygons}
   \end{figure}

   An \emph{integral unimodular affine transformation} of $\R^2$ is an affine map
   \begin{displaymath}
     \phi:\R^2\longrightarrow\R^2:\alpha\mapsto A\cdot \alpha+\tau
   \end{displaymath}
   with  $\tau\in\Z^2$ and $A\in\Gl_2(\Z)$ invertible over the integers.
   Such an integral unimodular affine
   transformation ${\phi}$ maps each face of $Q$ to a face of the convex lattice polygon ${\phi}(Q)$ and
   preserves thereby the number of lattice points on each
   face. Moreover, ${\phi}$ induces an isomorphism of the polarized toric
   surfaces $X_\ca$ and $X_{{\phi}(\ca)}$
   (cf.~\cite[Proposition~5.1.2]{GKZ94}). From the point of view of toric
   surfaces it therefore suffices to consider the marked
   polygon $(Q,\ca)$ only up to integral unimodular affine
   transformations, and if we suppose $\ca=Q\cap\Z^2$ then there are
   precisely sixteen of 
   them which we divide into two groups, $Q_a$, $Q_b$ and $Q_c$
   respectively $Q_{ca},\ldots,Q_{cm}$ (see Figure \ref{fig:polygons},
   cf.~\cite{Rab89} or \cite{PR00}). We fix the interior point at
   position $(1,1)$.

   The marked polygon $(Q_c,\ca_c)$ corresponds to $\P_\K^2$ embedded into
   $\P_\K^9$ via the $3$-uple Veronese embedding. The marked polygon $(Q_b,\ca_b)$
   corresponds to $\P_\K^1\times\P_\K^1$ embedded into $\P_\K^8$ via
   the $(2,2)$-Segre embedding. The marked polygon $(Q_a,\ca_a)$ describes the
   singular weighted projective plane $\P_\K(2,1,1)$ embedded into
   $\P^8$. 

   If a polygon marked polygon $(Q',\ca')$ is derived from $(Q,\ca)$ by cutting off one lattice
   point $(k,l)$, like $(Q_{ca},\ca_{ca})$ and $(Q_c,\ca_c)$, then the toric surface
   $X_{\ca'}$ is a blow up of $X_\ca$ 
   in a single point. Moreover, in 
   the torus $(\K^*)^2$ the hyperplane sections corresponding to
   \begin{displaymath}
     f=\sum_{(i,j)\in\ca'} a_{ij}\cdot x^iy^j=\sum_{(i,j)\in\ca}
       a_{ij}\cdot x^iy^j
   \end{displaymath}
   with $a_{kl}=0$ coincide. In particular, if they are both smooth
   their \emph{$j$-invariant} coincides since two birationally
   equivalent curves are already isomorphic
   (cf. \cite[Section~I.6]{Har77}). Since the $13$ polygons
   $Q_{ca},\ldots,Q_{cm}$ in 
   the second group in Figure \ref{fig:polygons} are all subpolygons
   of $Q_c$ the corresponding toric surfaces are all obtained from
   the projective plane by a couple of blow ups. When we want
   to compute the $j$-invariant of the curve corresponding to some
   Laurent polynomial with support in one of these $13$ polygons, we
   can instead consider the plane curve with support in $\ca_c$ but
   with the appropriate coefficients being zero. 

   Once we are able to compute the $j$-invariant for
   polynomials with support $\ca_a$, $\ca_b$ and $\ca_c$ we are 
   therefore able to compute the $j$-invariant for every Laurent
   polynomial with support on the lattice points of a lattice polygon
   with only one interior point. 

   We still assume that $(Q,\ca)$ is a marked lattice polygon with
   only one interior lattice point as above. Moreover, we use the notation
   $\underline{a}=(a_{ij}\;|\;(i,j)\in\ca)$, 
   and we suppose that
   \begin{displaymath}
     f=\sum_{(i,j)\in\ca}a_{ij}\cdot x^iy^j,
   \end{displaymath}
   then the $j$-invariant
   \begin{displaymath}
     j(C_f)=j(f)=\frac{A_\ca}{B_\ca}
   \end{displaymath}
   of the curve $C_f$ in $X_\ca$ defined by $f$
   can be expressed as a quotient of two homogeneous polynomials
   $A_\ca,B_\ca\in\Q[\underline{a}]$ of degree $12$. 
   In the case of $\ca=\ca_c$\; $A_\ca$ has $1607$ terms and $B_\ca$ has $2040$. In the
   case of $\ca=\ca_b$\; $A_\ca$ has $990$ terms and $B_\ca$ has
   $1010$. And finally in the case $\ca=\ca_a$\; $A_\ca$ has $267$ terms and
   $B_\ca$ has $312$. Every other case can be reduced to these
   three via some integral unimodular affine transformation and by
   setting some coefficients equal to zero.
   The reader interested in seeing or using the polynomials
   can consult the procedure \texttt{invariantsDB} in the Singular library
   \texttt{jinvariant.lib} (see \cite{KMM07}).
   The proof of our result relies heavily on the investigation of
   the combinatorics of these polynomials.

%%%%%%%%%%%%%%%%%%%%%%%%%%%%%%%%%%%%%%%%%%%%%%%%%%%%%%%%%%%%%%%%%%%%%%%%%

   \section{Tropicalization}\label{sec:tropicalization}

   In this section we want to pass from the algebraic to the tropical
   side. For this we specify a field $\K$ with a \emph{non-archimedean
     valuation} $\val:\K^*\rightarrow\R$ as base field and we extend
   the valuation to $\K$ by $\val(0)=\infty$.    
   We call $\val(k)$ also the \emph{tropicalization} of $k$.
   In the examples that we
   consider $\K$ will always be the \emph{field of Puiseux series}
   \begin{displaymath}
     \bigcup_{N=1}^\infty\Quot\Big(\C\big[\big[t^{\frac{1}{N}}\big]\big]\Big)
     =\left\{\sum_{\nu=m}^\infty c_\nu\cdot
       t^\frac{\nu}{N}\;\Big|\;
       c_\nu\in\C, N\in\Z_{>0}, m\in\Z\right\}
   \end{displaymath}
   and the valuation of a Puiseux series is its \emph{order}.

   If $f=\sum a_{ij}\cdot x^iy^j\in\K[x,y,x^{-1},y^{-1}]$ is
   any Laurent polynomial, we call the set
   \begin{displaymath}
     \supp(f)=\{(i,j)\in\Z^2\;|\;a_{ij}\not=0\}
   \end{displaymath}
   the \emph{support} of $f$ and the convex hull $\newton(f)$ of $\supp(f)$ in
   $\R^2$ is called the \emph{Newton polygon} of $f$. If
   $\supp(f)\subseteq \ca\subseteq\newton(f)\cap\Z^2$ then $f$ defines a curve $C_f$ in the
   toric surface $X_\ca$ as described in Section \ref{sec:toricsurfaces} and 
   we define the \emph{tropicalization} of $C_f$ as
   \begin{displaymath}
     \Trop\big(C_f\big)=\overline{\val\big(C_f\cap
       (\K^*)^2)}\subseteq\R^2,
   \end{displaymath}
   i.e.\  the closure of $\val\big(C_f\cap (\K^*)^2\big)$ with respect to
   the Euclidean topology in $\R^2$. Here by abuse of notation
   \begin{displaymath}
     \val:(\K^*)^2\longrightarrow\Q^2:(k_1,k_2)\mapsto
     \big(\val(k_1),\val(k_2)\big)
   \end{displaymath}
   denotes the Cartesian product of the above valuation map.

   A better way to compute the tropicalization of $C_f$ is as the
   tropical curve defined by the \emph{tropicalization} of the polynomial
   $f$, i.e.\ the piecewise linear map
   \begin{displaymath}
     \trop(f):\R^2\longrightarrow\R:(x,y)\mapsto\min\{\val(a_{ij})+i\cdot x+j\cdot y\;|\;(i,j)\in\supp(f)\}. 
   \end{displaymath}
   Given any \emph{plane tropical Laurent polynomial}
   \begin{displaymath}
     F:\R^2\longrightarrow\R:(x,y)\mapsto\min\{u_{ij}+i\cdot x+j\cdot y\;|\;(i,j)\in\ca'\}     
   \end{displaymath}
   with \emph{support} $\supp(F)=\ca'\subset\Z^2$ finite
   and $u_{ij}\in\R$, we call the \emph{locus $\mathcal{C}_F$ of
     non-differentiability} of $F$, i.e.\ the set of points 
   $(x,y)\in\R^2$ where the minimum is attained at least twice, the
   \emph{plane tropical curve} defined by $F$. The convex hull
   $\newton(F)$ of $\supp(F)$ is again called the \emph{Newton polygon} of $F$.

   By \emph{Kapranov's Theorem} (see \cite[Theorem~2.1.1]{EKL06}),
   $\Trop(C_f)$ coincides with the plane tropical curve defined by the
   plane tropical polynomial $\trop(f)$. In particular, $\Trop(C_f)$ is a
   piece-wise linear graph. 

   The plane tropical Laurent polynomial $F$ induces a \emph{marked subdivision}
   (cf.~\cite[Definition~7.2.1]{GKZ94}) of the
   marked polygon $\big(\newton(F),\ca\big)$ with
   $\supp(F)\subseteq\ca\subseteq\newton(F)\cap\Z^2$ in the following way: project the lower
   faces of the convex hull of
   \begin{displaymath}
     \{(i,j,u_{ij}\;|\;(i,j)\in\supp(F)\}
   \end{displaymath}
   into the $xy$-plane to subdivide of $\newton(F)$ into smaller
   polygons and mark those lattice points for which $(i,j,u_{ij})$ is
   contained in a lower face. 

   This subdivision is \emph{dual} to the
   tropical curve $\mathcal{C}_F$ in the following sense (see
   \cite[Prop.\ 3.11]{Mik05}): 
   Each marked polygon of the subdivision is dual to a
   vertex of $\mathcal{C}_F$, and each facet of a marked polygon is dual
   to an edge of $\mathcal{C}_F$. Moreover, if the facet, say $e$, has end points
   $(x_1,y_1)$ and $(x_2,y_2)$ then the \emph{direction vector} $v(E)$
   \label{page:directionvector} of the
   dual edge $E$ in $\mathcal{C}_F$ is defined (up to sign) as
   \begin{displaymath}
     v(E)=(y_2-y_1,x_1-x_2)^t
   \end{displaymath}
   and points in the direction of $E$.
   In particular, the edge $E$ is orthogonal to its dual facet $e$.
   Finally, the edge $E$ is unbounded if and
   only if its dual facet $e$ is contained in a
   facet of $\newton(F)$.

   \begin{example}\label{ex:curve}
     Consider the polynomial
     \begin{displaymath}
       f=xy+t\cdot (y+x+x^2+x^2y^2)+t^3
     \end{displaymath}
     The following diagram shows the support of $f$ and its marked Newton
     polygon.
     
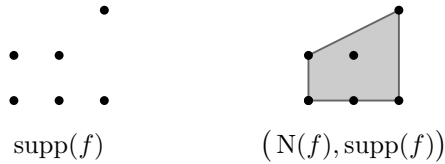
\begin{center}
       \begin{tabular}{c@{\hspace*{2cm}}c}
         \begin{texdraw}
           \drawdim cm  \relunitscale 0.6            \fcir f:0 r:0.04
           \move (2 2) 
           \fcir f:0 r:0.1
           \move (2 0) 
           \fcir f:0 r:0.1
           \move (1 1) 
           \fcir f:0 r:0.1
           \move (0 1) 
           \fcir f:0 r:0.1
           \move (1 0) 
           \fcir f:0 r:0.1
           \move (0 0) 
           \fcir f:0 r:0.1
         \end{texdraw}
         &
         \begin{texdraw}
           \drawdim cm  \relunitscale 0.6 \setgray 0.4
           \move (0 0)        
           \lvec (2 0)
           \lvec (2 2)
           \lvec (0 1)        
           \lvec (0 0)
           \lfill f:0.8
           \fcir f:0 r:0.1
           \move (2 2) 
           \fcir f:0 r:0.1
           \move (2 0) 
           \fcir f:0 r:0.1
           \move (1 1) 
           \fcir f:0 r:0.1
           \move (0 1) 
           \fcir f:0 r:0.1
           \move (1 0) 
           \fcir f:0 r:0.1
           \move (0 0) 
           \fcir f:0 r:0.1
%            \move (2 1) 
%            \fcir f:0 r:0.1
         \end{texdraw}
         \\[0.2cm]
         $\supp(f)$ & $\big(\newton(f),\supp(f)\big)$
       \end{tabular}       
     \end{center}
     The tropicalization of $f$ is
     \begin{displaymath}
       \trop(f):\R^2\rightarrow\R:(x,y)\mapsto\min\{x+y,1+y,1+x,1+2x,1+2x+2y,3\}.
     \end{displaymath}
     The support and Newton polygon of $f$ respectively of $\trop(f)$
     coincide. In order to compute the marked subdivision of the
     Newton polygon note that the points
     \begin{displaymath}
       (0,1,1),(1,0,1),(2,0,1),(2,2,1)
     \end{displaymath}
     lie in a plane while watching from below the point $(1,1,0)$ sticks out from this
     plane and the point $(0,0,3)$ lies way above it. We therefore get
     the following subdivision of the Newton polygon:
     
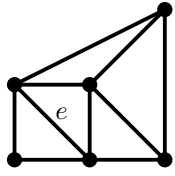
\begin{center}
       \begin{texdraw}
         \drawdim cm  \relunitscale 1 
         \linewd 0.05
         \move (2 2)        
         \lvec (2 0)
         \move (2 0)        
         \lvec (0 0)
         \move (0 0)        
         \lvec (0 1)
         \move (0 1)        
         \lvec (2 2)          
         \move (1 1)        
         \lvec (2 2)
         \move (2 0)        
         \lvec (1 1)
         \move (1 1)        
         \lvec (0 1)
         \move (1 0)        
         \lvec (1 1)
         \move (1 0)        
         \lvec (0 1)
         \move (2 2) 
         \fcir f:0 r:0.1
         \move (2 0) 
         \fcir f:0 r:0.1
         \move (1 1) 
         \fcir f:0 r:0.1
         \move (0 1) 
         \fcir f:0 r:0.1
         \move (1 0) 
         \fcir f:0 r:0.1
         \move (0 0) 
         \fcir f:0 r:0.1
         \htext (0.55 0.55){$e$}
       \end{texdraw}         
     \end{center}
     The polygon spanned by $(0,0)$, $(1,0)$ and $(0,1)$ is dual to
     the vertex of the tropical curve where the terms $3$, $1+x$ and
     $1+y$ take their common minimum, which is at the point
     $(x,y)=(2,2)$. Similarly the polygon spanned by $(1,0)$, $(1,1)$
     and $(0,1)$ corresponds to the point $(x,y)=(1,1)$, and the
     common face $e$ of the two polygons then is dual to the edge
     connecting these two points. Note that the direction vector of
     this edge $E$ is $v(E)=(1,1)$ is orthogonal to the
     face $e$ connecting the points $(1,0)$ and $(0,1)$ and points
     from the starting point $(1,1)$ of $E$ to its end point $(2,2)$. Computing the
     remaining vertices and edges of $\Trop(C_f)$ we get the following
     graph.
     
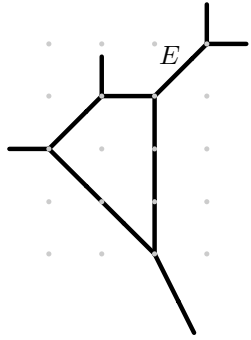
\begin{center}
       \begin{texdraw}
         \drawdim cm  \relunitscale 0.7 \arrowheadtype t:V
         % \linewd 0.05 \lpatt (0.1 0.4)
         % \move (-4 0) \avec (9 0) \move (0 -4) \avec (0 9)
         \linewd 0.08  \lpatt (1 0)
         %\setgray 0.6
         \relunitscale 1
         \move (-1 0) \fcir f:0 r:0.05
         \move (-1 0) \lvec (1 -2)
         \move (-1 0) \lvec (0 1)
         \move (-1 0) \rlvec (-0.75 0)
         \move (1 -2) \fcir f:0 r:0.05
         \move (1 -2) \lvec (1 1)
         \move (1 -2) \rlvec (0.75 -1.5)
         \move (0 1) \fcir f:0 r:0.05
         \move (0 1) \lvec (1 1)
         \move (0 1) \rlvec (0 0.75)
         \move (1 1) \fcir f:0 r:0.05
         \move (1 1) \lvec (2 2)
         \move (2 2) \fcir f:0 r:0.05
         \move (2 2) \rlvec (0 0.75)
         \move (2 2) \rlvec (0.75 0)
         \move (-1 -2) \fcir f:0.8 r:0.05
         \move (-1 -1) \fcir f:0.8 r:0.05
         \move (-1 0) \fcir f:0.8 r:0.05
         \move (-1 1) \fcir f:0.8 r:0.05
         \move (-1 2) \fcir f:0.8 r:0.05
         \move (0 -2) \fcir f:0.8 r:0.05
         \move (0 -1) \fcir f:0.8 r:0.05
         \move (0 0) \fcir f:0.8 r:0.05
         \move (0 1) \fcir f:0.8 r:0.05
         \move (0 2) \fcir f:0.8 r:0.05
         \move (1 -2) \fcir f:0.8 r:0.05
         \move (1 -1) \fcir f:0.8 r:0.05
         \move (1 0) \fcir f:0.8 r:0.05
         \move (1 1) \fcir f:0.8 r:0.05
         \move (1 2) \fcir f:0.8 r:0.05
         \move (2 -2) \fcir f:0.8 r:0.05
         \move (2 -1) \fcir f:0.8 r:0.05
         \move (2 0) \fcir f:0.8 r:0.05
         \move (2 1) \fcir f:0.8 r:0.05
         \move (2 2) \fcir f:0.8 r:0.05
         \htext (1.1 1.6){$E$}
       \end{texdraw}       
     \end{center}
   \end{example}
   
%%%%%%%%%%%%%%%%%%%%%%%%%%%%%%%%%%%%%%%%%%%%%%%%%%%%%%%%%%%%%%%%%%%%%%%%%

   \section{The tropical $j$-invariant of an elliptic plane tropical curve}\label{sec:jinv}

   For the purpose of this paper we want to define an elliptic plane
   tropical curve in the following way.

   \begin{definition}
     An \emph{elliptic plane tropical curve} is a tropical curve
     $\mathcal{C}_F$ defined by a plane tropical Laurent polynomial $F$ whose Newton
     polygon has precisely one interior lattice point. 
   \end{definition}

   The plane tropical curve $\mathcal{C}_F$ in Example \ref{ex:curve} is elliptic in
   this sense. Moreover, the graph $\mathcal{C}_F$ has \emph{genus}
   one, where the genus of a graph is the number of independent cycles
   of the graph. Obviously a cycle in the graph corresponds to an
   interior lattice point of the subdivision being a vertex of at
   least three polygons in the subdivision of the Newton
   polygon. We want to make this more precise in the following
   definition. 

   \begin{definition}
     Let $\mathcal{C}$ be a plane tropical curve with marked Newton
     polygon $(Q,\ca)$ and with dual marked
     subdivision $\{(Q_i,\ca_i)\;|\;i=1,\ldots,l\}$. Suppose
     that $\tilde{\omega}\in\Int(Q)\cap\Z^2$ and that the
     $(Q_i,\ca_i)$ are ordered such that $\tilde{\omega}$ is a vertex of $Q_i$
     for $i=1,\ldots,k$ and it is not contained in $Q_i$ for
     $i=k+1,\ldots,l$ (see Figure \ref{fig:cycle}). 
     \begin{figure}[h]
       \centering
       \begin{texdraw}
         \drawdim cm  \relunitscale 0.4 \linewd 0.05  \lpatt (1 0) \setgray 0
         \move (0 0) \rlvec (2 3) \rlvec (3 0)
         \lpatt (0.1 0.4) \rlvec (1.5 -1.5) \lpatt(1 0)
         \move (0 0) \rlvec (1 -3) \rlvec (4 0) \rlvec (1 1)
         \lpatt (0.1 0.4) \rlvec (1.3 1.3) \lpatt(1 0)
         \move (3 0) \lvec (0 0) \fcir f:0 r:0.1 \htext (-1.5 -0.2) {$\omega_2$}
         \move (3 0) \lvec (0.5 -1.5) \fcir f:0 r:0.1  \htext (-4.2 -2) {$\omega_1=\omega_{k+1}$}
         \move (3 0) \lvec (3.5 -3) \fcir f:0 r:0.1 
         \htext (2.5 -5.7) {$\xymatrix@R0.2cm{\omega_k\ar@{=}[d]\\\omega_0}$}
         \move (3 0) \lvec (5 -3) \fcir f:0 r:0.1  \htext (5.5 -4) {$\omega_{k-1}$}
         \move (3 0) \lvec (2 3) \fcir f:0 r:0.1 \htext (2 3.3) {$\omega_3$}
         \move (3 0) \lvec (5 3) \fcir f:0 r:0.1 \htext (5 3.3) {$\omega_4$}
         \move (3 0) \fcir f:0 r:0.1 \htext (3.7 -0.2) {$\tilde{\omega}$}
         \htext (1.5 -2.5) {$Q_k$}
         \htext (0.5 -1) {$Q_1$}
         \htext (1.2 0.8) {$Q_2$}
         \htext (3 1.8) {$Q_3$}
         \move (0 0) \rlvec (-1 3) \lvec (2 3) \htext (-0.7 2) {$Q_{k+1}$}
         \move (-1 3) \rlvec (3 3) \rlvec (3 -3) \htext (1.5 4.2) {$Q_{k+2}$}
         \move (2 6) \rlvec (2 0) \lpatt (0.1 0.4) \rlvec (2 0) 
         \move (-1 3) \lpatt (1 0) \rlvec (-1 -1) \lpatt (0.1 0.4)
         \rlvec (-1.5 -1.5)
         \move (1 -3) \lpatt (1 0) \rlvec (-0.5 -1) \lpatt (0.1 0.4)
         \rlvec (-0.5 -1)
         \move (5 -3) \lpatt (1 0) \rlvec (0 -1) \lpatt (0.1 0.4)
         \rlvec (0 -1)        
       \end{texdraw}    
       \caption{Marked subdivision determining a cycle}
       \label{fig:cycle}
     \end{figure}
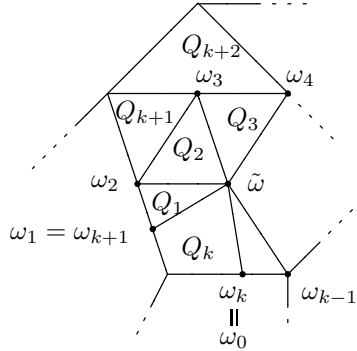
     We then say that $\tilde{\omega}$ \emph{determines a
       cycle} of $\mathcal{C}$, namely the union of the edges of $\mathcal{C}$ dual to the
     facets emanating from $\tilde{\omega}$, and we say that these edges 
     \emph{form the cycle} determined by $\tilde{\omega}$. We define
     the \emph{lattice length of the cycle} to be the sum of the
     lattice lengths of the edges which form the cycle, where for an
     edge $E$ with direction vector $v(E)$ (see
     p.~\pageref{page:directionvector}) the \emph{lattice length of $E$} is 
     \begin{displaymath}
       l(E)=\frac{||E||}{||v(E)||}
     \end{displaymath}
     the Euclidean length of $E$ divided by that of $v(E)$.
   \end{definition}

   \begin{example}\label{ex:curvecont}
     Coming back to our Example \ref{ex:curve} the curve has one cycle
     dual to the interior lattice point $(1,1)$ and it consists of
     four edges $E_1,\ldots,E_4$.
     
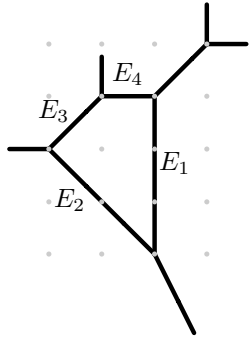
\begin{center}
       \begin{texdraw}
         \drawdim cm  \relunitscale 0.7 \arrowheadtype t:V
         % \linewd 0.05 \lpatt (0.1 0.4)
         % \move (-4 0) \avec (9 0) \move (0 -4) \avec (0 9)
         \linewd 0.08  \lpatt (1 0)
         %\setgray 0.6
         \relunitscale 1
         \move (-1 0) \fcir f:0 r:0.05
         \move (-1 0) \lvec (1 -2)
         \move (-1 0) \lvec (0 1)
         \move (-1 0) \rlvec (-0.75 0)
         \move (1 -2) \fcir f:0 r:0.05
         \move (1 -2) \lvec (1 1)
         \move (1 -2) \rlvec (0.75 -1.5)
         \move (0 1) \fcir f:0 r:0.05
         \move (0 1) \lvec (1 1)
         \move (0 1) \rlvec (0 0.75)
         \move (1 1) \fcir f:0 r:0.05
         \move (1 1) \lvec (2 2)
         \move (2 2) \fcir f:0 r:0.05
         \move (2 2) \rlvec (0 0.75)
         \move (2 2) \rlvec (0.75 0)
         \move (-1 -2) \fcir f:0.8 r:0.05
         \move (-1 -1) \fcir f:0.8 r:0.05
         \move (-1 0) \fcir f:0.8 r:0.05
         \move (-1 1) \fcir f:0.8 r:0.05
         \move (-1 2) \fcir f:0.8 r:0.05
         \move (0 -2) \fcir f:0.8 r:0.05
         \move (0 -1) \fcir f:0.8 r:0.05
         \move (0 0) \fcir f:0.8 r:0.05
         \move (0 1) \fcir f:0.8 r:0.05
         \move (0 2) \fcir f:0.8 r:0.05
         \move (1 -2) \fcir f:0.8 r:0.05
         \move (1 -1) \fcir f:0.8 r:0.05
         \move (1 0) \fcir f:0.8 r:0.05
         \move (1 1) \fcir f:0.8 r:0.05
         \move (1 2) \fcir f:0.8 r:0.05
         \move (2 -2) \fcir f:0.8 r:0.05
         \move (2 -1) \fcir f:0.8 r:0.05
         \move (2 0) \fcir f:0.8 r:0.05
         \move (2 1) \fcir f:0.8 r:0.05
         \move (2 2) \fcir f:0.8 r:0.05
         \htext (0.2 1.2){$E_4$}
         \htext (-1.2 0.5){$E_3$}
         \htext (-0.9 -1.2){$E_2$}
         \htext (1.1 -0.5){$E_1$}
       \end{texdraw}       
     \end{center}
     The edge $E_1$ is dual to the edge $e_1$ from $(0,1)$ to $(1,1)$
     in the Newton subdivision in Example \ref{ex:curve}, so that its
     direction vector is $v(E_1)=(0,-1)$ of Euclidean length $1$ and
     that the lattice length of $E_1$ is $l(E_1)=||E_1||=3$. Doing
     similar computations for the other edges the cycle length is
     \begin{displaymath}
       l(E_1)+l(E_2)+l(E_3)+l(E_4)=3+2+1+1=7.
     \end{displaymath}
   \end{example}
   
   \begin{definition}
     If $\mathcal{C}$ is an elliptic plane tropical curve then
     $\mathcal{C}$ has at most one cycle, and we define
     its \emph{tropical $j$-invariant} $j_{\trop}(\mathcal{C})$ to be
     the lattice length of this cycle if it has one. If $\mathcal{C}$
     has no cycle we define its tropical $j$-invariant to be zero.
   \end{definition}

   In Example \ref{ex:curvecont} the elliptic plane tropical curve has
   tropical $j$-invariant $7$.

   If a we fix the part of a Newton subdivision which determines the
   cycle then there is a nice formula to compute the cycle length, and
   thus the tropical $j$-invariant. For the proof we refer to
   \cite[Lemma~3.7]{KMM08}. 

  \begin{lemma}\label{lem:cyclelengthlinear}
    Let $(Q,\mathcal{A})$  be a marked polygon in $\R^2$
    with a marked subdivision $\{(Q_i,\ca_i)\;|\;i=1,\ldots,l\}$
    and suppose that $\tilde{\omega}\in\Int(Q)\cap\Z^2$ 
    is a vertex of $Q_i$ for $i=1,\ldots,k$ and it is not contained in
    $Q_i$ for $i=k+1,\ldots,l$. 
    
    If $u\in\R^\ca$ is such that the plane tropical curve
    \begin{displaymath}
      F=\min\{u_{ij}+i\cdot x+j\cdot y\;|\;(i,j)\in\mathcal{A}\}
    \end{displaymath}
    induces this subdivision (as described in Section
    \ref{sec:tropicalization}),  then $\tilde{\omega}$ determines a
    cycle in the plane  tropical curve $\mathcal{C}_F$ 
    and, using the notation in Figure \ref{fig:cycle}, its length is
    \begin{displaymath}
      \sum_{j=1}^k(u_{\tilde{\omega}}-u_{\omega_j})\cdot
      \frac{D_{j-1,j}+D_{j,j+1}+D_{j+1,j-1}}{D_{j-1,j}\cdot D_{j,j+1}}
    \end{displaymath}
    where $D_{i,j}=\det(w_i, w_j)$ with $w_i=\omega_i-\tilde{\omega}$ and $w_j=\omega_j-\tilde{\omega}$.
  \end{lemma}

  This formula implies in particular the following corollary.

  \begin{corollary}\label{cor:jinv}
    If $(Q,\ca)$ is a marked lattice polygon in $\R^2$ with precisely one interior
    lattice point, then
    \begin{displaymath}
      j_{\trop}:\R^\ca\longrightarrow\R:u\mapsto j_{\trop}(u):=j_{\trop}(\mathcal{C}_{F_u})
    \end{displaymath}
    with
    \begin{displaymath}
      F_u=\min\{u_{ij}+i\cdot x+j\cdot y\;|\;(i,j)\in\ca\}
    \end{displaymath}
    is a piecewise linear function which is linear on cones of
    the secondary fan (cf.~\cite[Chapter~7]{GKZ94}) of $\ca$. 
  \end{corollary}

%%%%%%%%%%%%%%%%%%%%%%%%%%%%%%%%%%%%%%%%%%%%%%%%%%%%%%%%%%%%%%%%%%%%%%%%%%%%%%%%%%%%%%%
   \section{Unimodular transformations preserve lattice length}\label{sec:unimodular}

   We want to relate the classical $j$-invariant to the
   tropical $j$-invariant, and we would again like to reduce the
   consideration of all possible Newton polygons with one interior
   point to the $16$ polygons in Figure \ref{fig:polygons}, or even
   better, to the three basic ones in the first group there. For that
   we have to understand the impact of an integral unimodular affine
   transformation on a plane tropical Laurent polynomial respectively
   the induced plane tropical curve.

   Given a linear form $l=u+i\cdot x+j\cdot y=u+(i,j)\cdot(x,y)^t$ with
   $i,j\in\Z$ and $u\in\R$ and given an integral unimodular affine
   transformation
   \begin{displaymath}
     \phi:\R^2\longrightarrow\R^2:\alpha\mapsto A\cdot \alpha+\tau
   \end{displaymath}
   with $A\in\Gl_2(\Z)$ and $\tau\in\Z^2$, we let $\phi$ act on $l$ via
   \begin{displaymath}
     l^\phi=u+(x,y)\cdot \phi\big((i,j)^t\big)
   \end{displaymath}
   and we let $\phi$ act on a plane tropical Laurent polynomial
   $F=\min\{u_{ij}+i\cdot x+j\cdot y\;|\;(i,j)\in\ca'\}$ via the
   linear forms, i.e.\
   \begin{displaymath}
     F^\phi=\min\{u_{ij}+(x,y)\cdot \phi\big((i,j)^t\big)\;|\;(i,j)\in\ca'\}.
   \end{displaymath}
   Note, that the translation by $\tau$ does not change
   the piecewise linear function defined by $F$ at all and the Newton
   polygon of $F$ is just translated by $\tau$. So $\tau$ has neither
   any impact on the Newton subdivision of $F$ nor on the tropical
   curve defined by $F$. Moreover, it is obvious that if
   $\{(Q_i,\ca_i)\;|\;i=1,\ldots,k\}$ is the marked subdivision of
   $\big(\newton(F),\supp(F))$ induced by $F$, then
   $\{\phi(Q_i),\phi(\ca_i)\;|\;i=1,\ldots,k\}$ is the marked
   subdivision of $\big(\newton(F^\phi),\supp(F^\phi)\big)$
   induced by $F^\phi$. 

   \kurz{
     It is a well-known fact that an integral unimodular affine
     transformation preserves lattice length, which implies the
     following corollary.
   }

   \begin{corollary}\label{cor:tropcurvetrans}
     Let $F$ be a plane tropical Laurent polynomial such that
     $\mathcal{C}_F$ is elliptic with positive tropical $j$-invariant
     and let $\phi$ be an integral unimodular affine transformation of
     $\R^2$, then $\mathcal{C}_{F^\phi}$ is elliptic with the same
     tropical $j$-invariant
     \begin{displaymath}
       j_{\trop}(\mathcal{C}_F)=j_{\trop}(\mathcal{C}_{F^\phi}).
     \end{displaymath}
   \end{corollary}

   \begin{example}
     Consider the polynomial
     \begin{displaymath}
       f=x^2y+xy^2+\frac{1}{t}\cdot xy+x+y
     \end{displaymath}
     inducing the following subdivision of its Newton polygon and the corresponding
     tropical curve:
     
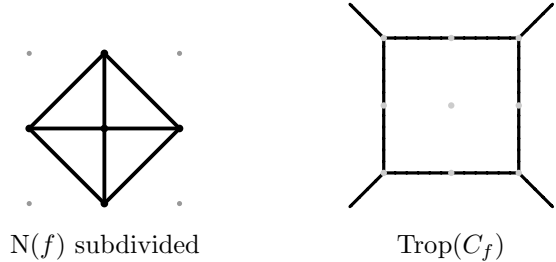
\begin{center}
       \begin{tabular}[m]{c@{\hspace*{2cm}}c}
         \begin{texdraw}
           \drawdim cm  \relunitscale 1
           \linewd 0.05
           \move (1 2)        
           \lvec (2 1)
           \move (2 1)        
           \lvec (1 0)
           \move (1 0)        
           \lvec (0 1)
           \move (0 1)        
           \lvec (1 2)      
           \move (2 1)        
           \lvec (1 1)
           \move (1 1)        
           \lvec (1 2)
           \move (1 0)        
           \lvec (1 1)
           \move (1 1)        
           \lvec (0 1)
           \move (0 0) \fcir f:0.6 r:0.03
           \move (0 1) \fcir f:0 r:0.05
           \move (0 2) \fcir f:0.6 r:0.03
           \move (1 0) \fcir f:0 r:0.05
           \move (1 1) \fcir f:0 r:0.05
           \move (1 2) \fcir f:0 r:0.05
           \move (2 0) \fcir f:0.6 r:0.03
           \move (2 1) \fcir f:0 r:0.05
           \move (2 2) \fcir f:0.6 r:0.03
         \end{texdraw}
         &
         \begin{texdraw}
           \drawdim cm  \relunitscale 0.15 \arrowheadtype t:V
           \linewd 0.3  \lpatt (1 0)       
           %\setgray 0.6
           \relunitscale 6
           \move (-1 -1) \fcir f:0 r:0.03
           \move (-1 -1) \lvec (-1 1)
           \move (-1 -1) \lvec (1 -1)
           \move (-1 -1) \rlvec (-0.5 -0.5)
           \move (-1 1) \fcir f:0 r:0.03
           \move (-1 1) \lvec (1 1)
           \move (-1 1) \rlvec (-0.5 0.5)
           \move (1 -1) \fcir f:0 r:0.03
           \move (1 -1) \lvec (1 1)
           \move (1 -1) \rlvec (0.5 -0.5)
           \move (1 1) \fcir f:0 r:0.03
           \move (1 1) \rlvec (0.5 0.5)         
           %% HERE STARTS THE CODE FOR THE LATTICE
           \move (-1 -1) \fcir f:0.8 r:0.05
           \move (-1 0) \fcir f:0.8 r:0.05
           \move (-1 1) \fcir f:0.8 r:0.05
           \move (0 -1) \fcir f:0.8 r:0.05
           \move (0 0) \fcir f:0.8 r:0.05
           \move (0 1) \fcir f:0.8 r:0.05
           \move (1 -1) \fcir f:0.8 r:0.05
           \move (1 0) \fcir f:0.8 r:0.05
           \move (1 1) \fcir f:0.8 r:0.05
           %% HERE ENDS THE CODE FOR THE LATTICE                          
         \end{texdraw}  
         \\[0.2cm]
         $\newton(f)$ subdivided & $\Trop(C_f)$
       \end{tabular}       
     \end{center}
     The plane tropical curve $\Trop(C_f)$ is elliptic with tropical
     $j$-invariant $8$. If we now
     apply the integral unimodular affine transformation
     \begin{displaymath}
       \phi:\R^2\longrightarrow \R^2:\alpha\mapsto 
       \begin{pmatrix}
         2 & 1\\1 &1 
       \end{pmatrix}
       \cdot \alpha
     \end{displaymath}
     we get
     \begin{displaymath}
       f^\phi=x^5y^3+x^4y^3+\frac{1}{t}\cdot x^3y^2+x^2y+xy
     \end{displaymath}
     with the following subdivision of its Newton polygon and
     the corresponding elliptic plane tropical curve having again
     tropical $j$-invariant $8$.
     
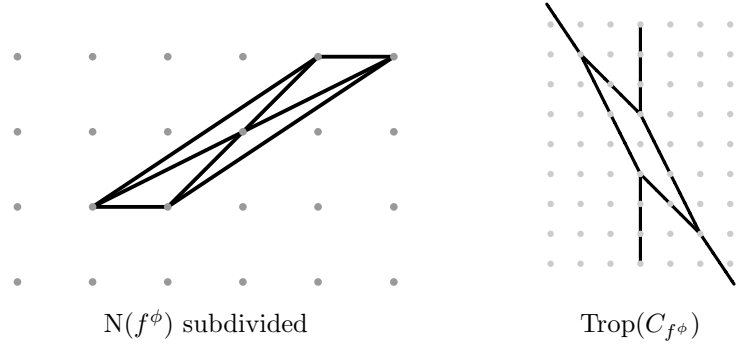
\begin{center}
       \begin{tabular}[m]{c@{\hspace*{2cm}}c}
         \begin{texdraw}
           \drawdim cm  \relunitscale 1
           \linewd 0.05
           \move (4 3)        
           \lvec (5 3)
           \move (5 3)        
           \lvec (2 1)
           \move (2 1)        
           \lvec (1 1)
           \move (1 1)        
           \lvec (4 3)    
           \move (5 3)        
           \lvec (3 2)
           \move (3 2)        
           \lvec (4 3)
           \move (2 1)        
           \lvec (3 2)
           \move (3 2)        
           \lvec (1 1)
           \move (0 0) \fcir f:0.6 r:0.05
           \move (0 1) \fcir f:0.6 r:0.05
           \move (0 2) \fcir f:0.6 r:0.05
           \move (0 3) \fcir f:0.6 r:0.05
           \move (1 0) \fcir f:0.6 r:0.05
           \move (1 1) \fcir f:0.6 r:0.05
           \move (1 2) \fcir f:0.6 r:0.05
           \move (1 3) \fcir f:0.6 r:0.05
           \move (2 0) \fcir f:0.6 r:0.05
           \move (2 1) \fcir f:0.6 r:0.05
           \move (2 2) \fcir f:0.6 r:0.05
           \move (2 3) \fcir f:0.6 r:0.05
           \move (3 0) \fcir f:0.6 r:0.05
           \move (3 1) \fcir f:0.6 r:0.05
           \move (3 2) \fcir f:0.6 r:0.05
           \move (3 3) \fcir f:0.6 r:0.05
           \move (4 0) \fcir f:0.6 r:0.05
           \move (4 1) \fcir f:0.6 r:0.05
           \move (4 2) \fcir f:0.6 r:0.05
           \move (4 3) \fcir f:0.6 r:0.05
           \move (5 0) \fcir f:0.6 r:0.05
           \move (5 1) \fcir f:0.6 r:0.05
           \move (5 2) \fcir f:0.6 r:0.05
           \move (5 3) \fcir f:0.6 r:0.05
         \end{texdraw}
         &
         \begin{texdraw}
           \drawdim cm  \relunitscale 0.15 \arrowheadtype t:V
           \linewd 0.3  \lpatt (1 0)
           %\setgray 0.6
           \relunitscale 2.66
           \move (0 -1) \fcir f:0 r:0.07
           \move (0 -1) \lvec (-2 3)
           \move (0 -1) \lvec (2 -3)
           \move (0 -1) \rlvec (0 -3)
           \move (-2 3) \fcir f:0 r:0.07
           \move (-2 3) \lvec (0 1)
           \move (-2 3) \rlvec (-1.125 1.685)
           \move (2 -3) \fcir f:0 r:0.07
           \move (2 -3) \lvec (0 1)
           \move (2 -3) \rlvec (1.125 -1.685)
           \move (0 1) \fcir f:0 r:0.07
           \move (0 1) \rlvec (0 3)           
           %% HERE STARTS THE CODE FOR THE LATTICE
           \move (-3 -4) \fcir f:0.8 r:0.1
           \move (-3 -3) \fcir f:0.8 r:0.1
           \move (-3 -2) \fcir f:0.8 r:0.1
           \move (-3 -1) \fcir f:0.8 r:0.1
           \move (-3 0) \fcir f:0.8 r:0.1
           \move (-3 1) \fcir f:0.8 r:0.1
           \move (-3 2) \fcir f:0.8 r:0.1
           \move (-3 3) \fcir f:0.8 r:0.1
           \move (-3 4) \fcir f:0.8 r:0.1
           \move (-2 -4) \fcir f:0.8 r:0.1
           \move (-2 -3) \fcir f:0.8 r:0.1
           \move (-2 -2) \fcir f:0.8 r:0.1
           \move (-2 -1) \fcir f:0.8 r:0.1
           \move (-2 0) \fcir f:0.8 r:0.1
           \move (-2 1) \fcir f:0.8 r:0.1
           \move (-2 2) \fcir f:0.8 r:0.1
           \move (-2 3) \fcir f:0.8 r:0.1
           \move (-2 4) \fcir f:0.8 r:0.1
           \move (-1 -4) \fcir f:0.8 r:0.1
           \move (-1 -3) \fcir f:0.8 r:0.1
           \move (-1 -2) \fcir f:0.8 r:0.1
           \move (-1 -1) \fcir f:0.8 r:0.1
           \move (-1 0) \fcir f:0.8 r:0.1
           \move (-1 1) \fcir f:0.8 r:0.1
           \move (-1 2) \fcir f:0.8 r:0.1
           \move (-1 3) \fcir f:0.8 r:0.1
           \move (-1 4) \fcir f:0.8 r:0.1
           \move (0 -4) \fcir f:0.8 r:0.1
           \move (0 -3) \fcir f:0.8 r:0.1
           \move (0 -2) \fcir f:0.8 r:0.1
           \move (0 -1) \fcir f:0.8 r:0.1
           \move (0 0) \fcir f:0.8 r:0.1
           \move (0 1) \fcir f:0.8 r:0.1
           \move (0 2) \fcir f:0.8 r:0.1
           \move (0 3) \fcir f:0.8 r:0.1
           \move (0 4) \fcir f:0.8 r:0.1
           \move (1 -4) \fcir f:0.8 r:0.1
           \move (1 -3) \fcir f:0.8 r:0.1
           \move (1 -2) \fcir f:0.8 r:0.1
           \move (1 -1) \fcir f:0.8 r:0.1
           \move (1 0) \fcir f:0.8 r:0.1
           \move (1 1) \fcir f:0.8 r:0.1
           \move (1 2) \fcir f:0.8 r:0.1
           \move (1 3) \fcir f:0.8 r:0.1
           \move (1 4) \fcir f:0.8 r:0.1
           \move (2 -4) \fcir f:0.8 r:0.1
           \move (2 -3) \fcir f:0.8 r:0.1
           \move (2 -2) \fcir f:0.8 r:0.1
           \move (2 -1) \fcir f:0.8 r:0.1
           \move (2 0) \fcir f:0.8 r:0.1
           \move (2 1) \fcir f:0.8 r:0.1
           \move (2 2) \fcir f:0.8 r:0.1
           \move (2 3) \fcir f:0.8 r:0.1
           \move (2 4) \fcir f:0.8 r:0.1
           \move (3 -4) \fcir f:0.8 r:0.1
           \move (3 -3) \fcir f:0.8 r:0.1
           \move (3 -2) \fcir f:0.8 r:0.1
           \move (3 -1) \fcir f:0.8 r:0.1
           \move (3 0) \fcir f:0.8 r:0.1
           \move (3 1) \fcir f:0.8 r:0.1
           \move (3 2) \fcir f:0.8 r:0.1
           \move (3 3) \fcir f:0.8 r:0.1
           \move (3 4) \fcir f:0.8 r:0.1
         \end{texdraw}
         \\[0.2cm]
         $\newton(f^\phi)$ subdivided & $\Trop(C_{f^\phi})$
       \end{tabular}
     \end{center}
   \end{example}

   The considerations in Section \ref{sec:toricsurfaces} together with
   this corollary will allow us to reduce the study of the
   tropicalization of an elliptic curve in a toric surface with an
   arbitrary Newton polygon with one interior point to the study of
   those whose Newton polygons are among the $16$ polygons in Figure
   \ref{fig:polygons}. 

%%%%%%%%%%%%%%%%%%%%%%%%%%%%%%%%%%%%%%%%%%%%%%%%%%%%%%%%%%%%%%%%%%%%%%%%%

   \section{The main result}\label{sec:main}

   Let us suppose now that $(Q,\ca)$ is a lattice polygon with only one
   interior point.

   \begin{remark}
     In Section \ref{sec:toricsurfaces} we have seen that 
     the $j$-invariant of a curve $C_f$
     with $\supp(f)\subseteq\ca$ can be computed by plugging the
     coefficients $a_{ij}$ of $f$ into a suitable quotient
     $j=\frac{A_\ca}{B_\ca}$ of homogeneous polynomials
     $A_\ca,B_\ca\in\Q[\underline{a}]$. This means in particular, that
     the valuation of the $j$-invariant can be read off $A_\ca$ and
     $B_\ca$ directly, unless some unlucky cancellation of leading
     terms occurs. 
   \end{remark}

   This leads to the following definition. 

   \begin{definition}
     The \emph{generic valuation} of a polynomial
     $0\not=H=\sum_\omega
     H_\omega\cdot\underline{a}^\omega\in\Q[\underline{a}]$ with
     $\underline{a}=(a_{ij}\;|\;(i,j)\in\ca)$ is 
     \begin{displaymath}
       \val_H:\R^{\ca}\longrightarrow\R:u\mapsto\val_H(u)=\min\{u\cdot \omega\;|\;H_\omega\not=0\},
     \end{displaymath}
     where
     \begin{displaymath}
       u\cdot \omega=\sum_{(i,j)\in\ca_c}u_{ij}\cdot \omega_{ij}.
     \end{displaymath}
     The \emph{generic valuation of the $j$-invariant} is the function
     \begin{displaymath}
       \val_j:\R^{\ca}\longrightarrow\R:u\mapsto\val_j(u)=\val_{A_\ca}(u)-\val_{B_\ca}(u).       
     \end{displaymath}
   \end{definition}

   Note that the tropical $j$-invariant is a \emph{tropical rational
     function} in the sense of \cite[Sec.~2.2]{Mik06} and
   \cite[Def.~3.1]{AR07}.

   \begin{remark}\label{rem:genericval}
     As mentioned above, unless some unlucky cancellation of the
     leading terms occurs for any $f=\sum_{(i,j)\in\ca}a_{ij}\cdot x^iy^j\in\K[x,y]$ with 
     $u_{ij}=\val(a_{ij})$ for all $(i,j)\in\ca$ we have
     \begin{displaymath}
       \val_j(u)=\val\big(j(f)\big).
     \end{displaymath}
     Note also, that if $D$ is a cone of the Gr\"obner fan of
     $A_\ca$ and $D'$ is a  cone of the Gr\"obner fan
     of $B_\ca$ then
     \begin{displaymath}
       {\val_j}_|:D\cap D'\longrightarrow \R
     \end{displaymath}
     is \emph{linear} by definition, and if both are top-dimensional,
     then no unlucky cancellation of leading terms can occur. In
     particular, the \emph{generic valuation of the $j$-invariant} $\val_j$
     is a piece linear function.
   \end{remark}

   We can now state the main result of our paper whose proof is
   discussed in the subsequent sections.

   \begin{theorem}\label{thm:main}
     Let $(Q,\ca)$ be a lattice polygon with only one interior lattice
     point.
     
     If $u\in\R^\ca$ is such that $\mathcal{C}_{F}$ with
     \begin{displaymath}
       F:\R^2\longrightarrow \R:(x,y)\mapsto\min\{u_{ij}+i\cdot
       x+j\cdot y\;|\;(i,j)\in\ca\}
     \end{displaymath}
     has a cycle, then
     \begin{displaymath}
       \val_j(u)=-j_{\trop}(u).
     \end{displaymath}
     Moreover, if $u$ is in a top-dimensional cone of the secondary
     fan of $\ca$ and $f=\sum_{(i,j)\in\ca}a_{ij}\cdot x^iy^j$ with
     $\val(a_{ij})=u_{ij}$, then
     \begin{displaymath}
       \val\big(j(f)\big)=-j_{\trop}(\mathcal{C}_F)=-j_{\trop}\big(\Trop(C_f)\big).
     \end{displaymath}
   \end{theorem}

   \begin{example}
     Consider the curve $C_f$ defined by
     \begin{displaymath}
       f=t^{\frac{3}{2}}\cdot(y+x^2+xy^2)+xy
     \end{displaymath}
     with the following subdivision of the Newton polygon and
     the corresponding elliptic plane tropical curve $\Trop(C_f)$:
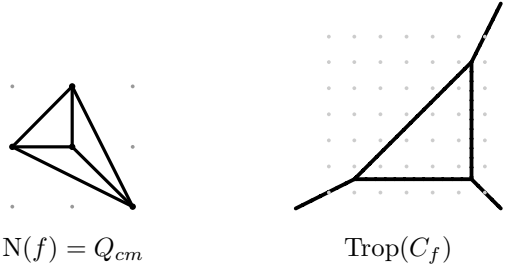
     \begin{center}
       \begin{tabular}[m]{c@{\hspace*{2cm}}c}
         \begin{texdraw}
           \drawdim cm  \relunitscale 0.8 
           \linewd 0.05
           \move (2 0)        
           \lvec (0 1)
           \move (0 1)        
           \lvec (1 2)
           \move (1 2)        
           \lvec (2 0)    
           \move (0 1)        
           \lvec (1 1)
           \move (1 1)        
           \lvec (2 0)
           \move (1 2)        
           \lvec (1 1)
           \move (0 0) \fcir f:0.6 r:0.03
           \move (0 1) \fcir f:0 r:0.05
           \move (0 2) \fcir f:0.6 r:0.03
           \move (1 0) \fcir f:0.6 r:0.03
           \move (1 1) \fcir f:0 r:0.05
           \move (1 2) \fcir f:0 r:0.05
           \move (2 0) \fcir f:0 r:0.05
           \move (2 1) \fcir f:0.6 r:0.03
           \move (2 2) \fcir f:0.6 r:0.03
           \move (2 0) 
           \fcir f:0 r:0.04
           \move (1 1) 
           \fcir f:0 r:0.04
           \move (0 1) 
           \fcir f:0 r:0.04
           \move (1 2) 
           \fcir f:0 r:0.04
         \end{texdraw}
         &
         \begin{texdraw}
           \drawdim cm  \relunitscale 0.13 \arrowheadtype t:V
           % \linewd 0.05 \lpatt (0.1 0.4)
           % \move (-4 0) \avec (9 0) \move (0 -4) \avec (0 9)
           \linewd 0.4  \lpatt (1 0)           
           %\setgray 0.6
           \relunitscale 2.66
           \move (1.5 3) \fcir f:0 r:0.07
           \move (1.5 3) \lvec (1.5 -1.5)
           \move (1.5 3) \lvec (-3 -1.5)
           \move (1.5 3) \rlvec (1.12 2.25)
           \move (1.5 -1.5) \fcir f:0 r:0.07
           \move (1.5 -1.5) \lvec (-3 -1.5)
           \move (1.5 -1.5) \rlvec (1.12 -1.12)
           \move (-3 -1.5) \fcir f:0 r:0.07
           \move (-3 -1.5) \rlvec (-2.25 -1.12)
           %% HERE STARTS THE CODE FOR THE LATTICE
           \move (-4 -2) \fcir f:0.8 r:0.07
           \move (-4 -1) \fcir f:0.8 r:0.07
           \move (-4 0) \fcir f:0.8 r:0.07
           \move (-4 1) \fcir f:0.8 r:0.07
           \move (-4 2) \fcir f:0.8 r:0.07
           \move (-4 3) \fcir f:0.8 r:0.07
           \move (-4 4) \fcir f:0.8 r:0.07
           \move (-3 -2) \fcir f:0.8 r:0.07
           \move (-3 -1) \fcir f:0.8 r:0.07
           \move (-3 0) \fcir f:0.8 r:0.07
           \move (-3 1) \fcir f:0.8 r:0.07
           \move (-3 2) \fcir f:0.8 r:0.07
           \move (-3 3) \fcir f:0.8 r:0.07
           \move (-3 4) \fcir f:0.8 r:0.07
           \move (-2 -2) \fcir f:0.8 r:0.07
           \move (-2 -1) \fcir f:0.8 r:0.07
           \move (-2 0) \fcir f:0.8 r:0.07
           \move (-2 1) \fcir f:0.8 r:0.07
           \move (-2 2) \fcir f:0.8 r:0.07
           \move (-2 3) \fcir f:0.8 r:0.07
           \move (-2 4) \fcir f:0.8 r:0.07
           \move (-1 -2) \fcir f:0.8 r:0.07
           \move (-1 -1) \fcir f:0.8 r:0.07
           \move (-1 0) \fcir f:0.8 r:0.07
           \move (-1 1) \fcir f:0.8 r:0.07
           \move (-1 2) \fcir f:0.8 r:0.07
           \move (-1 3) \fcir f:0.8 r:0.07
           \move (-1 4) \fcir f:0.8 r:0.07
           \move (0 -2) \fcir f:0.8 r:0.07
           \move (0 -1) \fcir f:0.8 r:0.07
           \move (0 0) \fcir f:0.8 r:0.07
           \move (0 1) \fcir f:0.8 r:0.07
           \move (0 2) \fcir f:0.8 r:0.07
           \move (0 3) \fcir f:0.8 r:0.07
           \move (0 4) \fcir f:0.8 r:0.07
           \move (1 -2) \fcir f:0.8 r:0.07
           \move (1 -1) \fcir f:0.8 r:0.07
           \move (1 0) \fcir f:0.8 r:0.07
           \move (1 1) \fcir f:0.8 r:0.07
           \move (1 2) \fcir f:0.8 r:0.07
           \move (1 3) \fcir f:0.8 r:0.07
           \move (1 4) \fcir f:0.8 r:0.07
           \move (2 -2) \fcir f:0.8 r:0.07
           \move (2 -1) \fcir f:0.8 r:0.07
           \move (2 0) \fcir f:0.8 r:0.07
           \move (2 1) \fcir f:0.8 r:0.07
           \move (2 2) \fcir f:0.8 r:0.07
           \move (2 3) \fcir f:0.8 r:0.07
           \move (2 4) \fcir f:0.8 r:0.07
           %% HERE ENDS THE CODE FOR THE LATTICE           
         \end{texdraw}
         \\[0.2cm]
         $\newton(f)=Q_{cm}$ 
         &
         $\Trop(C_f)$
       \end{tabular}
     \end{center}     
     The vertices of $\Trop(C_f)$ are
     \begin{displaymath}
       \left(\frac{3}{2},3\right),\;\;\;
       \left(\frac{3}{2},-\frac{3}{2}\right)\;\;\;\mbox{ and }\;\;\;\left(-3,-\frac{3}{2}\right),
     \end{displaymath}
     so that its
     tropical $j$-invariant is
     $j_{\trop}\big(\Trop(f)\big)=\frac{27}{2}$, while its $j$-invariant
     \begin{displaymath}
       j(f)=\frac{1+72\cdot t^{\frac{9}{2}}+1728\cdot
         t^9+13824\cdot t^{\frac{27}{2}}}{t^{\frac{27}{2}}+27\cdot t^{18}}
     \end{displaymath}
     has valuation $-\frac{27}{2}$.
   \end{example}

   An immediate consequence of the above theorem is the following corollary.

   \begin{corollary}
     If $f=\sum_{(i,j)\in\ca}a_{ij}\cdot x^iy^j\in\K[x,y]$ defines a
     smooth elliptic curve in $X_\ca$ whose $j$-invariant has
     non-negative valuation, then $\Trop(C_f)$ has no cycle.
   \end{corollary}

%%%%%%%%%%%%%%%%%%%%%%%%%%%%%%%%%%%%%%%%%%%%%%%%%%%%%%%%%%%%%%%%%%%%%%%%%%%%%%%%%%

   \section{Reduction to $\ca\in\{\ca_a,\ca_b,\ca_c\}$}\label{sec:reduction}

   Using an integral unimodular transformation we may assume that $Q$
   is one of the $16$ polygons in Figure \ref{fig:polygons}, since
   the application of such a transformation does not effect the
   statement of Theorem \ref{thm:main} due to Corollary
   \ref{cor:tropcurvetrans} and Section \ref{sec:toricsurfaces}.
   
   Next we want to reduce to the cases $\ca\in\{\ca_a,\ca_b,\ca_c\}$.

   If $f=\sum_{(i,j)\in\ca} a_{ij}\cdot x^iy^j\in\K[x,y]$ with $\supp(f)\subseteq\ca$ and we replace $f$
   by
   \begin{displaymath}
     f'=f+t^\alpha\cdot \sum_{(i,j)\in\ca\setminus\supp(f)}x^iy^j
   \end{displaymath}
   where $\alpha$ is much larger than
   $\max\{\val(a_{ij})\;|\;(i,j)\in\supp(f)\}$, then obviously
   \begin{displaymath}
     \val\big(j(f)\big)=\val\big(j(f')\big).
   \end{displaymath}

   Moreover, if we allow to plug in into $\val_j$ points $u$ where
   some of the $u_{ij}$ are $\infty$ (as long as the result still is
   a well defined real number), then we can evaluate $\val_j$ at $u$
   with $u_{ij}=\val(a_{ij})\in\R\cup\{\infty\}$ , $f$ defines a smooth elliptic
   curve and we get obviously
   \begin{displaymath}
     \val_j(u)=\val_j(u')
   \end{displaymath}
   where $u'_{ij}=u_{ij}$ for $(i,j)\in\supp(f)$ and else $u'_{ij}=\alpha$  with $\alpha$
   sufficiently large. 

   Finally, if in the definition of $F_u$ we allow some $u_{ij}$ to be
   $\infty$ then with the above notation the cycle of 
   $\mathcal{C}_{F_u}$ and $\mathcal{C}_{F_{u'}}$ will not change, so
   that
   \begin{displaymath}
     j_{\trop}(u)=j_{\trop}(\mathcal{C}_{F_u})=j_{\trop}(\mathcal{C}_{F_{u'}})=j_{\trop}(u'). 
   \end{displaymath}
   
   This shows that whenever we may as well assume that
   $\ca\in\{\ca_a,\ca_b,\ca_c\}$. 

%%%%%%%%%%%%%%%%%%%%%%%%%%%%%%%%%%%%%%%%%%%%%%%%%%%%%%%%%%%%%%%%%%%%%%%%%

   \section{The cases $\ca_a$, $\ca_b$ and $\ca_c$}\label{sec:cases}

   The case $\ca_c$ has been treated in \cite{KMM08}, and the two
   other cases work along the same lines. We therefore will be rather
   short in our presentation. Instead of considering all the cases by hand,
   as was done in \cite{KMM08} we will refer to computations done
   using the computer algebra systems \texttt{polymake} \cite{GJ97},
   \texttt{TOPCOM} \cite{Ram02} and \textsc{Singular} \cite{GPS05}. The
   code that we used for this is contained in the \textsc{Singular}  
   library \texttt{jinvariant.lib} (see \cite{KMM07}) and it is  
   available via the URL 
   \begin{center}
     http://www.mathematik.uni-kl.de/\textasciitilde keilen/en/jinvariant.html.
   \end{center}

   Fix now $(Q,\ca)$ with $\ca\in\{\ca_a,\ca_b,\ca_c\}$.

   We first of all observe that by Corollary \ref{cor:jinv} the
   tropical $j$-invariant is linear on the cones of the secondary fan
   of $\ca$ and that by Lemma \ref{lem:cyclelengthlinear} we can read
   off the
   assignment rule on each cone from the Newton subdivision of
   $(Q,\ca)$. Moreover, for the statement in Theorem \ref{thm:main} we
   only have to consider such cones for which the interior lattice
   point of $Q$ is visible in the subdivision. 

   If $U_\ca\subseteq\R^\ca$ is the union of
   these cones, then it was in each of the cases
   $\ca\in\{\ca_a,\ca_b,\ca_c\}$ computed by the procedure
   \texttt{testInteriorInequalities} 
   in the library \texttt{jinvariant.lib} that $U_\ca$ is contained in a
   single cone of the Gr\"obner fan of $A_\ca$ and  that
   \begin{displaymath}
     {\val_{A_\ca}}_|:U_\ca\longrightarrow\R:u\mapsto 12\cdot u_{11}.
   \end{displaymath}
   It suffices therefore to show that $\val_{B_\ca}$ is linear on
   the cones of the secondary fan of $\ca$ and to compare the
   assignment rules for $\val_j$ and $j_{\trop}$ on each of
   these cones.

   The two marked polygons $(Q_b,\ca_b)$ and $(Q_c,\ca_c)$ define \emph{smooth}
   toric surfaces and in these cases $B_\ca=\Delta_\ca$ is the
   $\ca$-discriminant of $\ca$ (cf.~\cite[Chapter~9]{GKZ94}). Therefore, by the
   Prime Factorization Theorem (see \cite[Theorem~10.1.2]{GKZ94})
   the secondary fan of $\ca_b$
   respectively $\ca_c$ is a refinement of the Gr\"obner fan of the
   $B_{\ca_b}$ respectively $B_{\ca_c}$.
   In view of Remark \ref{rem:genericval} and by the above
   considerations this shows in particular that $\val_j$ is linear on
   each cone of the secondary fan of $\ca$ for $\ca\in\{\ca_b,\ca_c\}$
   which is contained in $U_\ca$. The comparison of the assignment
   rules of the two linear functions $\val_j$ and $j_{\trop}$ on each
   of the cones contained in $U_\ca$ was done by the procedure
   \texttt{displayFan} from the \textsc{Singular} library
   \texttt{jinvariant.lib} using \texttt{TOPCOM} and
   \texttt{polymake}. It produces two postscript files which show all
   the different cases together with the assignment rules. The files
   are available via
   \begin{center}
     \small http://www.mathematik.uni-kl.de/\textasciitilde keilen/download/Tropical/secondary\_fan\_of\_2x2.ps
   \end{center}
   for $\ca=\ca_b$ respectively via
   \begin{center}
     \small http://www.mathematik.uni-kl.de/\textasciitilde keilen/download/Tropical/secondary\_fan\_of\_cubic.ps
   \end{center}
   for $\ca=\ca_c$ respectively. $849$ cases have to be considered for
   $\ca=\ca_c$ and $255$ for $\ca=\ca_b$.

   In the case $\ca=\ca_a$ the toric surface $X_\ca$ is \emph{not
     smooth}, but a quadric cone. Moreover, in this case
   $B_\ca$ is \emph{not} the $\ca$-discriminant $\Delta_\ca$, but instead
   \begin{displaymath}
     B_\ca=u_{02}^2\cdot \Delta_\ca.%=u_{02}\cdot D_\ca,
   \end{displaymath}
   Thus, the Gr\"obner fan of $B_\ca$ coincides with the Gr\"obner fan
   of $\Delta_\ca$ and it is still true by the Prime Factorization
   Theorem that the secondary fan of $\ca$ is a refinement of the
   Gr\"obner fan of $B_\ca$. We can therefore argue as above, and the
   case distinction ($202$ cases) can be viewed via
   \begin{center}
     \small http://www.mathematik.uni-kl.de/\textasciitilde keilen/download/Tropical/secondary\_fan\_of\_4x2.ps.
   \end{center}

   This finishes our proof, where for the ``moreover'' part we take
   Remark \ref{rem:genericval} into account.

   \begin{remark}
     It follows from the proof that $j_{\trop}$ is indeed linear on
     each cone of the Gr\"obner fan of $B_\ca$ in the above cases.
     This could have been proved directly
     with the same argument as in \cite[Lemma~5.2]{KMM08}. 

     If we denote by $D_\ca$ the regular $\ca$-determinant 
     (cf.~\cite[Section~11.1]{GKZ94}), then $D_{\ca_b}=\Delta_{\ca_b}$
     and $D_{\ca_c}=\Delta_{\ca_c}$ by \cite[Theorem~11.1.3]{GKZ94}
     since $X_\ca$ is smooth in these cases. Even though in general
     the regular $\ca$-determinant is not a polynomial, it is so for
     $\ca=\ca_a$ by \cite[Theorem~11.1.6]{GKZ94} since $\ca_a$ is
     quasi-smooth by \cite[Theorem~5.4.12]{GKZ94} in the sense of that
     theorem. More precisely, we have
     \begin{displaymath}
       D_{\ca_a}=u_{02}\cdot \Delta_{\ca_a}
     \end{displaymath}
     and by \cite[Theorem~11.1.3]{GKZ94} it is a divisor of the
     principal $\ca$-determinant $E_\ca$
     (cf.~\cite[Chapter~9]{GKZ94}). Therefore, the secondary fan of
     $\ca$ (which is 
     the Gr\"obner fan of $E_\ca$) is a refinement of the Gr\"obner fan of
     $D_\ca$ and thus of the Gr\"obner fan of $B_\ca$. 

     One therefore could have used the description of the vertices of the Newton
     polytope of $D_\ca$ in \cite[Theorem~11.3.2]{GKZ94} in order to
     show that on each cone of the Gr\"obner fan of $B_\ca$ contained
     in $U_\ca$ the two functions $\val_j$ and $j_{\trop}$ coincide by
     a direct case study as was done in \cite[Lemma~5.5]{KMM08} for
     the case $\ca=\ca_c$.
   \end{remark}

%%%%%%%%%%%%%%%%%%%%%%%%%%%%%%%%%%%%%%%%%%%%%%%%%%%%%%%%%%%%%%%%%%%%%%%%%

%  \bibliographystyle{plain}
%  \bibliography{bibliographie}

\begin{thebibliography}{10}

\bibitem{AR07}
Lars Allermann and Johannes Rau.
\newblock First steps in tropical intersection theory.
\newblock arXiv:0709.3705, 2007.

\bibitem{EKL06}
Manfred Einsiedler, Mikhail Kapranov, and Douglas Lind.
\newblock Non-archimedean amoebas and tropical varieties.
\newblock {\em J.\ Reine Angew.\ Math.}, 601:139--157, 2006.

\bibitem{Ful93}
William Fulton.
\newblock {\em Introduction to Toric Varieties}.
\newblock Princeton University Press, 1993.

\bibitem{GJ97}
Ewgenij Gawrilow and Michael Joswig.
\newblock {\sc polymake} 2.3.
\newblock Technical report, TU Berlin and TU Darmstadt, 1997.
\newblock http://www.math.tu-berlin.de/polymake.

\bibitem{GKZ94}
Israel~M. Gelfand, Mikhail~M. Kapranov, and Andrei~V. Zelevinsky.
\newblock {\em Discriminants, Resultants, and Multidimensional Determinants}.
\newblock Birkh\"auser, 1994.

\bibitem{GPS05}
G.-M. Greuel, G.~Pfister, and H.~Sch\"onemann.
\newblock {\sc Singular} 3.0.
\newblock {A Computer Algebra System for Polynomial Computations}, Centre for
  Computer Algebra, University of Kaiserslautern, 2005.
\newblock {\tt http://www.singular.uni-kl.de}.

\bibitem{Har77}
Robin Hartshorne.
\newblock {\em Algebraic Geometry}.
\newblock Springer, 1977.

\bibitem{JMM07a}
Anders~Nedergaard Jensen, Hannah Markwig, and Thomas Markwig.
\newblock tropical.lib.
\newblock A {\sc singular} 3.0 library for computations in tropical geometry,
  2007.
\newblock {\tt http://www.singular.uni-kl.de/\textasciitilde
  keilen/de/tropical.html}.

\bibitem{KMM08}
Eric Katz, Hannah Markwig, and Thomas Markwig.
\newblock The $j$-invariant of a plane tropical cubic.
\newblock Preprint, 2007.

\bibitem{KMM07}
Eric Katz, Hannah Markwig, and Thomas Markwig.
\newblock jinvariant.lib.
\newblock A {\sc singular} 3.0 library for computations with $j$-invariants in
  tropical geometry, 2007.
\newblock {\tt http://www.singular.uni-kl.de/\textasciitilde
  keilen/de/jinvariant.html}.

\bibitem{KM08}
Michael Kerber and Hannah Markwig.
\newblock Counting tropical elliptic curves with fixed $j$-invariant.
\newblock math.AG/0608472, 2006.

\bibitem{Mik05}
Grigory Mikhalkin.
\newblock Enumerative tropical geometry in {${\mathbbm R^2}$}.
\newblock {\em J. Amer. Math. Soc.}, 18:313--377, 2005.
\newblock math.AG/0312530.

\bibitem{Mik06}
Grigory Mikhalkin.
\newblock Tropical geometry and its applications.
\newblock In {\em International Congress of Mathematicians}, volume~II, pages
  827--852. Eur.\ Math.\ Soc., 2006.

\bibitem{PR00}
Bjorn Poonen and Fernando Rodriguez-Villegas.
\newblock Lattice polygons and the number 12.
\newblock {\em Amer. Math. Monthly}, 107(3):238--250, 2000.

\bibitem{Rab89}
Stanley Rabinowitz.
\newblock A census of convex lattice polygons with at most one interior lattice
  point.
\newblock {\em Ars Combin.}, 28:83--96, 1989.

\bibitem{Ram02}
J\"org Rambau.
\newblock Topcom: Triangulations of point configurations and oriented matroids.
\newblock In Xiao-Shan~Gao Arjeh~M.~Cohen and Nobuki Takayama, editors, {\em
  Mathematical Software - ICMS 2002}, pages 330--340, 2002.
\newblock available at
  http://www.uni-bayreuth.de/departments/wirtschaftsmathematik/rambau/TOPCOM/.

\bibitem{Spe07}
David Speyer.
\newblock Uniformizing tropical curves i: Genus zero and one.
\newblock arXiv:0711.2677, 2007.

\bibitem{Vig04}
Magnus~Dehli Vigeland.
\newblock The group law on a tropical elliptic curve.
\newblock math.AG/0411485, 2004.

\end{thebibliography}

\end{document}